\documentclass[11pt,a4paper]{article}

\usepackage[T1]{fontenc}
\usepackage[utf8]{inputenc}
\usepackage{lmodern}
\usepackage[margin=29mm]{geometry}
\usepackage{amsmath,amssymb,amsthm,mathtools}
\usepackage{microtype}
\usepackage{enumitem}
\usepackage[hidelinks]{hyperref}

\setlist[enumerate]{label=\textup{(\roman*)},leftmargin=*,itemsep=3pt,topsep=4pt}
\allowdisplaybreaks[1]
\numberwithin{equation}{section}

\newtheorem{theorem}{Theorem}[section]
\newtheorem{lemma}[theorem]{Lemma}
\newtheorem{proposition}[theorem]{Proposition}
\newtheorem{corollary}[theorem]{Corollary}
\theoremstyle{remark}

\newcommand{\Q}{Q}
\newcommand{\F}{\mathcal F}
\newcommand{\A}{\mathcal A}

\newcommand{\C}{\mathcal C}
\newcommand{\J}{\mathcal J}
\newcommand{\PP}{\mathcal P}
\newcommand{\SSS}{\mathcal S}
\newcommand{\NN}{\mathcal N}
\newcommand{\eps}{\varepsilon}
\DeclareMathOperator{\diam}{diam}

\DeclareMathOperator{\supp}{supp}
\newcommand{\doi}[1]{\href{https://doi.org/#1}{doi:\,\nolinkurl{#1}}}

\hypersetup{
  pdftitle={Binary Codes with a Forbidden Distance},
  pdfauthor={Jianwei Jiang and Chunhua Yang},
  pdfsubject={Forbidden Hamming distances and extremal binary codes},
  pdfkeywords={binary code, forbidden distance, Hamming cube, vertex isoperimetry}
}

\title{Binary Codes with a Forbidden Distance}
\author{Jianwei Jiang\qquad Chunhua Yang\textsuperscript{*}\\[5pt]
  \small School of Mathematics and Statistics, Weifang University,\\
  \small Weifang, Shandong, China\\[3pt]
  \small \textsuperscript{*}Corresponding author: \texttt{chunhuayangwf@sina.com}}
\date{}

\begin{document}
\maketitle

\begin{abstract}
The maximum size of a $t$-intersecting code in $[m]^n$ was determined independently by Ahlswede and Khachatrian and by Frankl and Tokushige,
resolving a problem posed by Frankl and F\"uredi. A code is $t$-intersecting if every two distinct words agree in at least $t$ coordinates. Keevash, Lifshitz, Long and Minzer studied the weaker $(t-1)$-avoiding condition and proved that, for $m>2$ and $n$ sufficiently large in terms of $t$, $(t-1)$-avoiding codes have the same maximum size and extremal ball constructions as $t$-intersecting codes. They left the binary case open.

We solve the binary case for all $t\ge3$ and $n\ge(t+1)^2$. If the forbidden Hamming distance $n-t+1$ is odd, the maximum size of a $(t-1)$-avoiding code strictly exceeds that of a $t$-intersecting code. If it is even, the two maxima coincide, but their extremal families differ: some maximum $(t-1)$-avoiding codes are not $t$-intersecting. We also characterize all maximum $(t-1)$-avoiding codes in the even case.

For odd forbidden distances, a weight-parity construction and a fixed-point-free pairing determine the maximum size, while Kleitman’s diameter theorem gives the comparison with $t$-intersecting codes. For even forbidden distances, parity splitting and coordinate deletion reduce the problem to avoiding two consecutive distances. A component decomposition and disjoint-neighbourhood arguments, together with Harper’s vertex-isoperimetric theorem, yield the sharp upper bound attained by a $t$-intersecting construction. Frankl’s equality characterization for the diameter theorem then identifies the extremal families, and elementary binomial estimates establish the stated range of $n$.
\end{abstract}


\section{Introduction}

Determining the maximum size of a $t$-intersecting code in $[m]^n$
is a classical problem in extremal combinatorics.
Here $[m]=\{1,\ldots,m\}$, and a code
$\F\subseteq[m]^n$ is called $t$-intersecting
if every two distinct words agree in at least $t$ coordinates.
This problem was posed by Frankl and
F\"uredi~\cite{FranklFuredi1980}
and solved independently by Ahlswede and
Khachatrian~\cite{AhlswedeKhachatrian1998}
and by Frankl and Tokushige~\cite{FranklTokushige1999},
who determined the maximum size and identified the
corresponding extremal constructions.

A natural weakening of the $t$-intersecting condition is to forbid only
pairs of words that agree in exactly $t-1$ coordinates. Such codes are
called $(t-1)$-avoiding. Keevash, Lifshitz, Long and
Minzer~\cite[Theorem~1.1]{KLLM2023} proved that, for $m\ge3$ and $n$
sufficiently large in terms of $t$, $(t-1)$-avoiding codes have the same maximum size as $t$-intersecting codes, and the same extremal ball constructions up to isomorphism. They also asked whether the same conclusion holds for $m=2$, but their junta approximation method does not apply to the conjectured binary extremal ball constructions~\cite[Section~9]{KLLM2023}. More recently, Iarovikova, Noskov, Sokolov and Terekhov~\cite[Theorem~3]{INST2026} established the exact extremal result under polynomial lower bounds on both $m$
and $n$ in terms of $t$, but their theorem does not cover $m=2$.

In this paper, we consider the binary case. Let $\Q_n=\{0,1\}^n$ be the
binary Hamming cube, and let $d_H(x,y)$ denote the number of coordinates
in which $x,y\in\Q_n$ differ. A binary code is any subset of $\Q_n$;
no linearity is assumed. Since $x$ and $y$ agree in exactly
$n-d_H(x,y)$ coordinates, the $(t-1)$-avoiding condition
is equivalent to forbidding Hamming distance $n-t+1$.

The distinction between $t$-intersecting and $(t-1)$-avoiding codes is
important. In the binary setting, the former condition requires the
code to have diameter at most $n-t$, whereas the latter excludes only
the single distance $n-t+1$, allowing pairs at larger distances.
Kleitman's diameter theorem~\cite{Kleitman1966} determines the maximum
size under the diameter restriction, but does not directly give an
upper bound for the more general avoiding problem. Thus, although every
$t$-intersecting code is $(t-1)$-avoiding, it is not immediate whether
the two conditions yield the same maximum size. Related problems
concerning prescribed or forbidden Hamming distances have been studied
by Frankl~\cite{Frankl1986}, Enomoto, Frankl, Ito and
Nomura~\cite{EnomotoEtAl1987}, and Bassalygo, Cohen and
Z\'emor~\cite{BassalygoEtAl2000}.

\subsection{Main results}

Let $M_t(n)$ denote the maximum size of a binary code in $\Q_n$ avoiding
Hamming distance $n-t+1$.
The answer depends on the parity of the forbidden distance.

\begin{theorem}\label{thm:main}
Let $t\ge3$ and $n\ge(t+1)^2$ be integers. Then
\begin{equation}\label{eq:main}
 M_t(n)=
 \begin{cases}
  2^{n-1},& n-t+1\text{ is odd},\\[4pt]
  \displaystyle 2\sum_{j=0}^{(n-t-1)/2}\binom{n-1}{j},
       & n-t+1\text{ is even}.
 \end{cases}
\end{equation}
\end{theorem}

When the forbidden distance is odd, the maximum size of a $(t-1)$-avoiding code
strictly exceeds the maximum size of a $t$-intersecting code. When it is even,
the two maxima coincide. However, as our next result shows, their
extremal families need not be the same.

For a positive integer $N$ and $z\in\Q_N$, let $|z|$ denote its Hamming
weight, and let $\pi(z)=|z|\bmod2\in\{0,1\}$. Write
$B_r(c)=\{z\in\Q_N:d_H(z,c)\le r\}$ for the Hamming ball with centre $c$
and radius $r$. With the last coordinate distinguished, define
\begin{equation}\label{eq:iota}
 \iota_\eps:\Q_N\longrightarrow\Q_{N+1},\qquad
 \iota_\eps(z)=(z,\eps\oplus\pi(z)),\qquad \eps\in\{0,1\},
\end{equation}
where $\oplus$ denotes addition modulo $2$, applied coordinatewise to binary
words of the same length. The image of $\iota_\eps$ is the set of words of weight congruent
to $\eps$ modulo $2$.

\begin{theorem}\label{thm:equality}
Let $t\ge3$ and $n\ge(t+1)^2$ be integers, and suppose that $n-t+1$ is even.
Put $N=n-1$ and $r=(N-t)/2$. A code $\F\subseteq\Q_n$ avoiding
distance $n-t+1$ has size $M_t(n)$ if and only if
\begin{equation}\label{eq:extremal}
 \F=\iota_0\bigl(B_r(c_0)\bigr)
       \cup\iota_1\bigl(B_r(c_1)\bigr)
 \qquad\text{for some }c_0,c_1\in\Q_N.
\end{equation}
There is no restriction relating $c_0$ and $c_1$.
\end{theorem}

If $c_0=c_1=c$, then $\F=B_r(c)\times\{0,1\}$ has diameter at most
$2r+1=n-t$ and is therefore $t$-intersecting. However, independent
choices of $c_0$ and $c_1$ can also produce maximum avoiding codes that
are not $t$-intersecting. Thus, unlike the $m\ge3$ case
in~\cite{KLLM2023}, the two problems have the same maximum size but
different extremal families when the forbidden distance is even. We
return to this distinction in Section~\ref{sec:threshold}.

\section{Preliminaries on the Hamming cube}\label{sec:prelim}

For an integer $N\ge1$, write $[N]=\{1,\ldots,N\}$. We identify a
word $x\in\Q_N$ with its support
$\supp(x)=\{i\in[N]:x_i=1\}\subseteq[N]$. Under this identification,
$d_H(x,y)=|\supp(x)\triangle\supp(y)|$, where $\triangle$ denotes
symmetric difference. We use the same notation for a family of words
and its corresponding family of supports.

For an integer $s$, put $V_N(s)=\sum_{j=0}^{s}\binom Nj$,
where an empty sum is $0$ and $\binom Nj=0$ for $j<0$ or $j>N$.
Thus $V_N(s)=0$ for $s<0$, $V_N(s)=2^N$ for $s\ge N$, and
$|B_s(c)|=V_N(s)$ for $s\ge0$.
For a family $X\subseteq\Q_N$ and an integer $h\ge0$, let
\[
 \NN_h(X)=\{z\in\Q_N:d_H(z,x)\le h\text{ for some }x\in X\}
\]
be its closed $h$-neighbourhood. In particular,
$\NN_h(\varnothing)=\varnothing$ and $\NN_0(X)=X$.
The ambient dimension will always be specified when this notation is
used. For a nonempty family $X$, write
$\diam(X)=\max\{d_H(x,y):x,y\in X\}$, and set
$\diam(\varnothing)=0$.

The antipode of $x\in\Q_N$ is
$\bar x=x\oplus\mathbf1$, where $\mathbf1$ is the all-one word;
write $\bar X=\{\bar x:x\in X\}$. We will use
\begin{equation}\label{eq:antipode}
 d_H(x,\bar y)=N-d_H(x,y),\qquad
 \NN_h(\bar X)=\overline{\NN_h(X)}.
\end{equation}
Here and throughout, $\bar X$ denotes the family obtained by complementing
each word of $X$.

We use the even-diameter case of Kleitman's diameter
theorem~\cite{Kleitman1966} (see also
Frankl~\cite[(3)--(4)]{Frankl2017}), together with Frankl's
equality characterization~\cite[Theorem~1]{Frankl2017}.

\begin{theorem}[Kleitman; Frankl]\label{thm:diameter}
Let $r\ge0$ and $N\ge2r+2$ be integers. If $X\subseteq\Q_N$
has diameter at most $2r$, then $|X|\le V_N(r)$. Equality holds if
and only if $X=B_r(c)$ for some $c\in\Q_N$.
\end{theorem}

The following lemma follows from the one-step neighbourhood estimate in
\cite[Lemma~4.1]{JiangYang2026}, whose proof uses Harper's
vertex-isoperimetric theorem~\cite{Harper1966}.

\begin{lemma}\label{lem:expansion}
Let $N\ge1$, $h\ge0$ and $0\le R\le N$ be integers. If
$X\subseteq\Q_N$ and $|X|\le V_N(R)$, then
\begin{equation}\label{eq:expansion}
 |\NN_h(X)|\ge\frac{V_N(R+h)}{V_N(R)}\,|X|.
\end{equation}
\end{lemma}

\begin{proof}
The assertion is immediate if $h=0$ or $R=N$.
We may therefore assume that $h\ge1$ and $R<N$.

We first establish a bound for closed $1$-neighbourhoods.
Let $s$ be an integer with $0\le s<N$, and let $Y\subseteq\Q_N$.
Choose $Y'\subseteq Y$ with $|Y'|=\min\{|Y|,V_N(s)\}$, and set
$A=Y'\times\{0,1\}\subseteq\Q_{N+1}$.
Since $|A|\le2V_N(s)$, applying
\cite[Lemma~4.1]{JiangYang2026} with $d=N+1$ and $r=s$ gives
$|\NN_1(A)|\ge V_N(s+1)|A|/V_N(s)$, where the neighbourhood of
$A$ is taken in $\Q_{N+1}$.
Since $\NN_1(A)=\NN_1(Y')\times\{0,1\}$
and $Y'\subseteq Y$, we obtain
\begin{equation}\label{eq:one-step}
 |\NN_1(Y)| \ge \frac{V_N(s+1)}{V_N(s)}\min\{|Y|,V_N(s)\}.
\end{equation}
Taking $Y=X$ and $s=R$ in \eqref{eq:one-step}, and using
the hypothesis $|X|\le V_N(R)$, proves \eqref{eq:expansion}
for $h=1$.

We now proceed by induction. Suppose that, for some integer $j\ge1$,
\[
 |\NN_j(X)| \ge \frac{V_N(R+j)}{V_N(R)}\,|X|.
\]
Put $X_j=\NN_j(X)$ and $s=R+j$.
If $s<N$, then $\NN_1(X_j)=\NN_{j+1}(X)$, so
\eqref{eq:one-step} gives
\[
 |\NN_{j+1}(X)| \ge
 \frac{V_N(s+1)}{V_N(s)}\min\{|X_j|,V_N(s)\}
 \ge \frac{V_N(s+1)}{V_N(R)}\,|X|.
\]
For the second inequality, the induction hypothesis and
$|X|\le V_N(R)$ show that both entries in the minimum are
at least $V_N(s)|X|/V_N(R)$.

If $s\ge N$, then $V_N(s)=V_N(s+1)=2^N$.
Since $X_j\subseteq\NN_{j+1}(X)$,
the induction hypothesis again gives the required bound.

This completes the induction and the proof.
\end{proof}

\section{Parity reduction and a lower bound}\label{sec:parity}

For $\eps\in\{0,1\}$, let
$E_\eps=\{x\in\Q_n:|x|\equiv\eps\pmod2\}$.
For $x,y\in\Q_n$, we have $d_H(x,y)\equiv |x|+|y|\pmod 2$. Hence two
words of the same weight parity are always at even Hamming distance.

\begin{proposition}\label{prop:odd}
If $1\le d\le n$ is odd, then the maximum size of a subset of
$\Q_n$ avoiding distance $d$ is $2^{n-1}$.
\end{proposition}

\begin{proof}
The sets $E_0$ and $E_1$ each have size $2^{n-1}$, and any two words
in the same set have the same weight parity, hence are at even Hamming
distance. Thus either $E_0$ or $E_1$ avoids the odd distance $d$, giving
the lower bound. Choose a word
$v\in\Q_n$ of weight $d$. Since $(x\oplus v)\oplus v=x$ and $v\ne0$,
the map $x\mapsto x\oplus v$ is a fixed-point-free involution. Moreover,
$d_H(x,x\oplus v)=|v|=d$. Hence its orbits partition $\Q_n$ into
$2^{n-1}$ pairs, each consisting of two words at distance $d$. A family
avoiding distance $d$ contains at most one word from each pair, giving
the upper bound $2^{n-1}$.
\end{proof}

We now consider a family $\F\subseteq\Q_n$ with a single forbidden
even distance, say $2r+2$. Parity splitting and coordinate deletion are also used in the
diameter problem; see~\cite[Section~4]{Frankl2017}.
The following lemma gives the precise reduction needed here.

\begin{lemma}\label{lem:projection}
Let $N=n-1$, let $r\in\mathbb{Z}_{\ge0}$, and suppose that $2r+2\le n$.
For $\F\subseteq\Q_n$, use the maps defined in \eqref{eq:iota} to define
$\A_\eps=\iota_\eps^{-1}(\F\cap E_\eps)$, $\eps\in\{0,1\}$.
Then $\F$ avoids distance $2r+2$ if and only if each of
$\A_0$ and $\A_1$ avoids both distances $2r+1$ and $2r+2$
in $\Q_N$. Moreover, $|\F|=|\A_0|+|\A_1|$.
\end{lemma}

\begin{proof}
Since $d_H(x,y)\equiv |x|+|y|\pmod 2$, any pair with one word
in $E_0$ and the other in $E_1$ is at odd
Hamming distance, and hence cannot have distance $2r+2$. It therefore
remains to consider pairs within a fixed parity class $E_\eps$.

Let $x,y\in\Q_N$, and write $s=d_H(x,y)$. If $s$ is even, then $x$
and $y$ have the same weight parity, so the last coordinates appended
by $\iota_\eps$ are equal. Hence
$d_H(\iota_\eps(x),\iota_\eps(y))=s$.
If $s$ is odd, then $x$ and $y$ have opposite weight parities, so the
appended last coordinates are different. Hence
$d_H(\iota_\eps(x),\iota_\eps(y))=s+1$.
Consequently,
\[
 d_H\bigl(\iota_\eps(x),\iota_\eps(y)\bigr)=2r+2
 \quad\Longleftrightarrow\quad
 s\in\{2r+1,2r+2\}.
\]
Thus $\F\cap E_\eps$ avoids distance $2r+2$ exactly when
$\A_\eps$ avoids both $2r+1$ and $2r+2$. Since each $\iota_\eps$ is a
bijection from $\Q_N$ onto $E_\eps$, we also have
$|\F|=|\F\cap E_0|+|\F\cap E_1|=|\A_0|+|\A_1|$.
\end{proof}

\begin{proposition}\label{prop:lower}
Let $t\ge3$ and $n\ge t+1$ be integers such that $n-t+1$ is even.
Put $N=n-1$ and $r=(N-t)/2$. For any $c_0,c_1\in\Q_N$,
the code $\F=\iota_0(B_r(c_0))\cup\iota_1(B_r(c_1))$ avoids distance
$n-t+1$ and has size $2V_N(r)$. Consequently, $M_t(n)\ge2V_N(r)$.
\end{proposition}

\begin{proof}
Set $\A_\eps=B_r(c_\eps)$ for $\eps\in\{0,1\}$. Each $\A_\eps$ has diameter at most $2r$ and hence avoids both $2r+1$ and $2r+2$. Since $2r+2=n-t+1$, Lemma~\ref{lem:projection} shows that
$\F=\iota_0(\A_0)\cup\iota_1(\A_1)$ avoids distance
$n-t+1$ in $\Q_n$. As the two lifted families lie in the disjoint classes $E_0$ and $E_1$
and each $\iota_\eps$ is a bijection,
$|\F|=|\A_0|+|\A_1|=2V_N(r)$. Thus $M_t(n)\ge2V_N(r)$.
\end{proof}

\section{Families avoiding two consecutive distances}\label{sec:components}

Throughout this section, let $t\ge3$, $r\ge0$ and $N=2r+t$
be integers. Suppose that $\A\subseteq\Q_N$ satisfies
\begin{equation}\label{eq:two-forbidden}
 d_H(x,y)\notin\{2r+1,2r+2\}
 \qquad\text{for all distinct }x,y\in\A.
\end{equation}
Define a graph $\Gamma_\A$ on $\A$ by joining distinct words
when their Hamming distance is $1$ or $2$. Let
$\C_1,\ldots,\C_s$ be its nonempty connected components;
when $\A=\varnothing$, the list is empty.

\subsection{The component decomposition}

\begin{lemma}\label{lem:gap}
Every component $\C_i$ has diameter at most $2r$. For distinct
components $\C_i,\C_j$, exactly one of the following holds:
\begin{enumerate}
 \item $d_H(x,y)\le2r$ for every $x\in\C_i$ and $y\in\C_j$;
 \item $d_H(x,y)\ge2r+3$ for every $x\in\C_i$ and $y\in\C_j$.
\end{enumerate}
\end{lemma}

\begin{proof}
Fix a component $\C_i$ and let $x,y\in\C_i$. Take a path
$x=y_0,y_1,\ldots,y_\ell=y$ in $\Gamma_\A$. If $d_H(x,y)>2r$,
let $j$ be the smallest index such that $d_H(x,y_j)>2r$. Then
$d_H(x,y_{j-1})\le2r$, while \eqref{eq:two-forbidden} gives
$d_H(x,y_j)\ge2r+3$. Hence
\[
 3\le d_H(x,y_j)-d_H(x,y_{j-1})
   \le d_H(y_{j-1},y_j)\le2,
\]
a contradiction. Thus $d_H(x,y)\le2r$, and since $x,y$ were
arbitrary, $\diam(\C_i)\le2r$.

Now let $\C_i$ and $\C_j$ be distinct components. If no pair
$x\in\C_i$, $y\in\C_j$ satisfies $d_H(x,y)\le2r$, then
\eqref{eq:two-forbidden} gives $d_H(x,y)\ge2r+3$ for every such
pair. Otherwise, choose $x_0\in\C_i$ and $y_0\in\C_j$ with
$d_H(x_0,y_0)\le2r$.

For any $x\in\C_i$, take a path
$x_0=x'_0,x'_1,\ldots,x'_p=x$ in $\C_i$. If
$d_H(x'_k,y_0)\le2r$, then
\[
 d_H(x'_{k+1},y_0)
 \le d_H(x'_k,y_0)+d_H(x'_k,x'_{k+1})
 \le2r+2.
\]
Since the distances $2r+1$ and $2r+2$ are forbidden, this implies
$d_H(x'_{k+1},y_0)\le2r$. Repeating this argument along the path gives
$d_H(x,y_0)\le2r$ for every $x\in\C_i$. Fixing such an $x$ and
applying the same argument along a path in $\C_j$ from $y_0$ to an
arbitrary $y\in\C_j$ yields $d_H(x,y)\le2r$. Thus the lemma follows.
\end{proof}

Let $\J\subseteq\{1,\ldots,s\}$ consist of the indices $i$ for
which alternative~\textup{(ii)} of Lemma~\ref{lem:gap} holds with
at least one other component. Define
\begin{equation}\label{eq:decomposition}
 \PP=\bigcup_{i\in\J}\C_i,\qquad
 \SSS=\A\setminus\PP,\qquad
 p=|\PP|,\quad q=|\SSS|,\quad a=V_N(r).
\end{equation}
Both $\PP$ and $\SSS$ are unions of whole components.
This decomposition is defined for every family satisfying
\eqref{eq:two-forbidden}.

\begin{lemma}\label{lem:small-components}
For each $i\in\J$, there is a word $z_i\in\Q_N$ such that
$\C_i\subseteq B_{t-3}(z_i)$. In particular, $|\C_i|\le V_N(t-3)$ for every $i\in\J$.
Moreover,
\begin{equation}\label{eq:cross-small}
 d_H(x,y)\le2r\qquad(x\in\SSS,\ y\in\A),
\end{equation}
and therefore $\diam(\SSS)\le2r$ and $q\le a$.
\end{lemma}

\begin{proof}
Fix $i\in\J$. By the definition of $\J$, there is a component $\C_j$
such that $d_H(x,y)\ge 2r+3$ for all $x\in\C_i$ and $y\in\C_j$.
Choose $y\in\C_j$ and set $z_i=\bar y$. Since $N=2r+t$, for every
$x\in\C_i$,
\[
 d_H(x,z_i)=d_H(x,\bar y)=N-d_H(x,y)
 \le N-(2r+3)=t-3.
\]
Hence $\C_i\subseteq B_{t-3}(z_i)$, and therefore
$|\C_i|\le V_N(t-3)$.

Now let $x\in\SSS$, and let $\C_i$ be the component containing $x$.
Since $i\notin\J$, the component $\C_i$ has no component at distance
at least $2r+3$ from it. Thus, by Lemma~\ref{lem:gap},
$d_H(x,y)\le2r$ for every $y$ in any component $\C_j\ne\C_i$.
The same bound holds for $y\in\C_i$, since
$\diam(\C_i)\le2r$. This proves \eqref{eq:cross-small}.
In particular, $\diam(\SSS)\le2r$. By the assumptions at the beginning of this section, $N=2r+t\ge2r+3$, so Theorem~\ref{thm:diameter} gives $q=|\SSS|\le V_N(r)=a$.
\end{proof}

\subsection{Disjoint neighbourhoods}

\begin{lemma}
The sets $\NN_1(\SSS)$, $\NN_{t-2}(\bar\SSS)$, and
$\NN_1(\C_i)$ for $i\in\J$ are pairwise disjoint. Consequently,
\begin{equation}\label{eq:packing}
 |\NN_1(\SSS)|+|\NN_{t-2}(\SSS)|
   +\sum_{i\in\J}|\NN_1(\C_i)|\le2^N.
\end{equation}
\end{lemma}

\begin{proof}
Suppose first that the closed one-neighbourhoods of two distinct
components $\C_i$ and $\C_j$ intersect. Then there are
$x\in\C_i$, $y\in\C_j$, and $z\in Q_N$ with
$d_H(x,z)\le1$ and $d_H(y,z)\le1$, so $d_H(x,y)\le2$.
Thus $x$ and $y$ are adjacent in $\Gamma_\A$, contradicting that
$\C_i$ and $\C_j$ are distinct components. Hence
$\NN_1(\SSS)$ and the sets $\NN_1(\C_i)$, $i\in\J$, are pairwise
disjoint.

It remains to separate $\NN_{t-2}(\bar\SSS)$. Suppose, for
contradiction, that $\NN_{t-2}(\bar\SSS)\cap\NN_1(\A)\ne\varnothing$.
Then there are $x\in\A$, $y\in\SSS$, and $z\in Q_N$ with
$d_H(z,x)\le1$ and $d_H(z,\bar y)\le t-2$, and hence
$d_H(x,\bar y)\le t-1$. On the other hand, \eqref{eq:cross-small}
and $N=2r+t$ give
$d_H(x,\bar y)=N-d_H(x,y)\ge t$, a contradiction. Since
$\NN_1(\SSS)$ and every $\NN_1(\C_i)$ are contained in
$\NN_1(\A)$, all the sets in the statement are pairwise disjoint.
Finally, antipodal complementation preserves neighbourhood size by
\eqref{eq:antipode}, so
$|\NN_{t-2}(\bar\SSS)|=|\NN_{t-2}(\SSS)|$.
Since these neighbourhoods are pairwise disjoint subsets of $Q_N$,
the sum of their cardinalities is at most $2^N$.
Together with the preceding identity, this gives \eqref{eq:packing}.
\end{proof}

\subsection{The upper bound}

Set
\begin{equation}\label{eq:kappa}
 \lambda_{N,t}=\frac{V_N(t-2)}{V_N(t-3)},\qquad
 \kappa_{N,t}=\frac{V_N(r)\lambda_{N,t}}{2^N}-1.
\end{equation}
The denominator $V_N(t-3)$ is positive because $t\ge3$.

\begin{proposition}\label{prop:quantitative}
Every family $\A\subseteq\Q_N$ satisfying
\eqref{eq:two-forbidden}, with the decomposition
\eqref{eq:decomposition}, satisfies
\begin{equation}\label{eq:quantitative}
 |\A|\le a-\kappa_{N,t}p.
\end{equation}
In particular, if $\kappa_{N,t}>0$, then $|\A|\le a$,
and equality implies $\PP=\varnothing$ and $\diam(\A)\le2r$.
\end{proposition}

\begin{proof}
For $i\in\J$, Lemma~\ref{lem:small-components} gives
$|\C_i|\le V_N(t-3)$. Apply Lemma~\ref{lem:expansion}
with $R=t-3$ and $h=1$ to obtain
\begin{equation}\label{eq:component-expansion}
 |\NN_1(\C_i)|\ge\lambda_{N,t}|\C_i|.
\end{equation}
Also, $q=|\SSS|\le a=V_N(r)$. Applying the same lemma to
$\SSS$ with $R=r$ and, respectively, $h=1$ and $h=t-2$ gives
\begin{equation}\label{eq:S-expansion}
 |\NN_1(\SSS)|+|\NN_{t-2}(\SSS)|
   \ge\frac{V_N(r+1)+V_N(r+t-2)}{a}\,q.
\end{equation}
Since $N=2r+t$, we have $r+t-2=N-r-2$. By binomial symmetry,
\[
 V_N(r+t-2)
 =\sum_{j=0}^{N-r-2}\binom Nj
 =\sum_{j=r+2}^{N}\binom Nj.
\]
Together with $V_N(r+1)=\sum_{j=0}^{r+1}\binom Nj$, this gives
\begin{equation}\label{eq:reflection}
 V_N(r+1)+V_N(r+t-2)=2^N.
\end{equation}

Substituting \eqref{eq:component-expansion} and
\eqref{eq:S-expansion} into \eqref{eq:packing}, and using
\eqref{eq:reflection} and $\sum_{i\in\J}|\C_i|=p$, yields
\[
 \frac{2^N}{a}\,q+\lambda_{N,t}p\le2^N.
\]
Multiplication by $a/2^N$ and addition of $p$ give
\[
 |\A|=q+p
 \le a-\left(\frac{a\lambda_{N,t}}{2^N}-1\right)p
 =a-\kappa_{N,t}p.
\]
If $\kappa_{N,t}>0$ and $|\A|=a$, then $p=0$.
Hence $\A=\SSS$, whose diameter is at most $2r$ by
Lemma~\ref{lem:small-components}.
\end{proof}

\begin{corollary}\label{cor:criterion}
Let $t\ge3$, $r\ge0$ and $N=2r+t$ be integers. If
\begin{equation}\label{eq:criterion}
 V_N(r)V_N(t-2)>2^NV_N(t-3),
\end{equation}
then the largest family in $\Q_N$ avoiding both $2r+1$ and
$2r+2$ has size $V_N(r)$. The extremal families are exactly
the balls $B_r(c)$, $c\in\Q_N$.
\end{corollary}

\begin{proof}
By the definitions in \eqref{eq:kappa},
$\kappa_{N,t}=V_N(r)V_N(t-2)/(2^N V_N(t-3))-1$,
so \eqref{eq:criterion} is equivalent to $\kappa_{N,t}>0$.
Proposition~\ref{prop:quantitative} therefore gives
$|\A|\le a-\kappa_{N,t}p\le a=V_N(r)$.
On the other hand, every ball $B_r(c)$ has diameter at most $2r$
and cardinality $V_N(r)$, so it is admissible and attains the upper bound.

If $|\A|=V_N(r)$, then $\kappa_{N,t}>0$ forces $p=0$ in
\eqref{eq:quantitative}. Thus $\A=\SSS$ and
$\diam(\A)\le2r$ by Lemma~\ref{lem:small-components}.
Since $N=2r+t\ge2r+3$, the equality case of
Theorem~\ref{thm:diameter} gives $\A=B_r(c)$ for some
$c\in\Q_N$. Hence the families attaining the maximum $V_N(r)$ are precisely the Hamming balls of radius $r$.
\end{proof}

\section{Dimension bound and extremal families}\label{sec:threshold}

We next prove that \eqref{eq:criterion} holds whenever
$N+1\ge(t+1)^2$. Only elementary estimates for binomial coefficients
are needed.

\begin{lemma}\label{lem:central}
For every integer $N\ge1$,
\begin{equation}\label{eq:central}
 \max_{0\le j\le N}\frac{\binom Nj}{2^N}
   \le\frac1{\sqrt{N+1}}.
\end{equation}
\end{lemma}

\begin{proof}
The binomial coefficients increase up to the middle layer and then decrease. Thus, if $N=2m$, the maximum is $\binom{2m}{m}$, while if $N=2m+1$, the two middle coefficients $\binom{2m+1}{m}$ and $\binom{2m+1}{m+1}$ are equal and maximal. We treat the two cases separately.

First let $N=2m$ and put $b_m=\binom{2m}{m}/4^m$. We prove by induction that $b_m\le(2m+1)^{-1/2}$. The claim is immediate for $m=0$. If it holds for $m$, then
\[
 b_{m+1}=\frac{2m+1}{2m+2}\,b_m
 \le\frac{\sqrt{2m+1}}{2m+2}
 \le\frac1{\sqrt{2m+3}},
\]
where the last inequality is equivalent to $(2m+1)(2m+3)\le(2m+2)^2$. Hence \eqref{eq:central} holds for even $N$.

Now let $N=2m+1$. Since the middle coefficients are maximal,
\[
 \max_{0\le j\le N}\frac{\binom Nj}{2^N}
 =\frac{\binom{2m+1}{m}}{2^{2m+1}}
 =\frac{\binom{2m+2}{m+1}}{2^{2m+2}}
 =b_{m+1}
 \le\frac1{\sqrt{2m+3}}
 \le\frac1{\sqrt{N+1}}.
\]
Thus \eqref{eq:central} also holds for odd $N$.
\end{proof}

\begin{lemma}\label{lem:growth}
If $N\ge1$ and $0\le k<N$ are integers, then
\begin{equation}\label{eq:growth}
 (N-k)V_N(k)\le(k+1)V_N(k+1).
\end{equation}
\end{lemma}

\begin{proof}
Count ordered pairs $(U,W)$ with $U,W\subseteq[N]$, $|U|\le k$,
$U\subset W$, and $|W\setminus U|=1$. By the definition of $V_N(k)$,
there are exactly $V_N(k)$ possible sets $U$. For each such $U$, there
are $N-|U|\ge N-k$ choices for $W$. Hence the number of these ordered
pairs is at least $(N-k)V_N(k)$.

Conversely, every such $W$ has size at most $k+1$. Thus there are at
most $V_N(k+1)$ possible sets $W$, and for each fixed $W$, at most
$k+1$ choices for $U$, obtained by deleting one element of $W$. Hence
the number of ordered pairs is at most $(k+1)V_N(k+1)$.
Comparing the two counts gives \eqref{eq:growth}.
\end{proof}

\begin{proposition}\label{prop:threshold}
Let $t\ge3$, $r\ge0$ and $N=2r+t$ be integers. If
$N+1\ge(t+1)^2$, then
\[
 \kappa_{N,t}\ge\frac{2t+5}{(t-2)(t+1)}>0.
\]
\end{proposition}

\begin{proof}
Recall from \eqref{eq:kappa} that
$\kappa_{N,t}=V_N(r)\lambda_{N,t}/2^N-1$, so we estimate the two
factors separately. Set $a=V_N(r)$. Since
$2^N=\sum_{j=0}^{N}\binom Nj$, we may decompose
\[
 2^N
 =\sum_{j=0}^{r}\binom Nj
 +\sum_{j=r+1}^{N-r-1}\binom Nj
 +\sum_{j=N-r}^{N}\binom Nj.
\]
By symmetry of the binomial coefficients, the first and third sums
are both equal to $a$, so
$2^N=2a+\sum_{j=r+1}^{N-r-1}\binom Nj$. Dividing by $2^{N+1}$ and
rearranging gives
$\frac{a}{2^N}=\frac12-\frac12\sum_{j=r+1}^{N-r-1}\frac{\binom Nj}{2^N}$.
Since $N=2r+t$, the middle sum contains exactly $t-1$ terms. By
Lemma~\ref{lem:central} and $N+1\ge(t+1)^2$,
\[
 \frac{a}{2^N}
 \ge\frac12-\frac{t-1}{2\sqrt{N+1}}
 \ge\frac1{t+1}.
\]
Next apply Lemma~\ref{lem:growth} with $k=t-3$. Since $0\le t-3<N$,
\[
 \lambda_{N,t}=\frac{V_N(t-2)}{V_N(t-3)}
 \ge\frac{N-t+3}{t-2}.
\]
Substituting these two bounds into \eqref{eq:kappa}, and using
$N+1\ge(t+1)^2$, we obtain
\[
 \kappa_{N,t}
 \ge\frac{N-t+3}{(t-2)(t+1)}-1
 \ge\frac{(t+1)^2-t+2}{(t-2)(t+1)}-1
 =\frac{2t+5}{(t-2)(t+1)}>0,
\]
since $t\ge3$.
\end{proof}

\begin{proof}[Proof of Theorem~\ref{thm:main}]
Let $d=n-t+1$. The assumptions give $1\le d\le n$.
If $d$ is odd, Proposition~\ref{prop:odd} gives
$M_t(n)=2^{n-1}$.

Suppose that $d$ is even. Put $N=n-1$ and
$r=(n-t-1)/2$, so that $r\ge0$, $N=2r+t$ and
$d=2r+2$. Given a code $\F\subseteq\Q_n$ avoiding $d$,
let $\A_0,\A_1\subseteq\Q_N$ be its two projections from
Lemma~\ref{lem:projection}. Each avoids $2r+1$ and $2r+2$.
Since $N+1=n\ge(t+1)^2$, Proposition~\ref{prop:threshold}
gives $\kappa_{N,t}>0$, and
Proposition~\ref{prop:quantitative} gives
$|\A_0|,|\A_1|\le V_N(r)$. Thus, by Lemma~\ref{lem:projection},
$|\F|=|\A_0|+|\A_1|\le2V_N(r)$.
Proposition~\ref{prop:lower} supplies a code of that size.
This proves the even case of \eqref{eq:main} and completes the proof.
\end{proof}

\begin{proof}[Proof of Theorem~\ref{thm:equality}]
Put $N=n-1$, $r=(n-t-1)/2$, and $a=V_N(r)$, so that
$N=2r+t$. Since $N+1=n\ge(t+1)^2$,
Proposition~\ref{prop:threshold} gives $\kappa_{N,t}>0$,
which is equivalent to \eqref{eq:criterion}.

Suppose that $\F\subseteq\Q_n$ avoids distance $n-t+1$ and has
size $|\F|=2a=M_t(n)$. By Lemma~\ref{lem:projection}, its
projections $\A_0,\A_1\subseteq\Q_N$ avoid both $2r+1$ and
$2r+2$, and $|\A_0|+|\A_1|=2a$.
Corollary~\ref{cor:criterion} gives $|\A_0|,|\A_1|\le a$,
so both families have size $a$. By the equality characterization
in that corollary, $\A_\eps=B_r(c_\eps)$ for some
$c_\eps\in\Q_N$, for each $\eps\in\{0,1\}$.
Lifting the two families yields \eqref{eq:extremal}.

Conversely, for any $c_0,c_1\in\Q_N$,
Proposition~\ref{prop:lower} shows that the code
in \eqref{eq:extremal} avoids distance $n-t+1$ and has
size $2a=M_t(n)$.
\end{proof}

We now compare the $(t-1)$-avoiding and $t$-intersecting problems under
the assumptions of Theorem~\ref{thm:main}.

Suppose first that $n-t+1$ is odd, and put $u=(n-t)/2$. A binary
$t$-intersecting code has diameter at most $2u$. Since
$n=2u+t\ge 2u+3$, Theorem~\ref{thm:diameter} bounds its size by
$V_n(u)$. Moreover,
\[
 2^n
 =\sum_{j=0}^{u}\binom nj
 +\sum_{j=u+1}^{n-u-1}\binom nj
 +\sum_{j=n-u}^{n}\binom nj.
\]
By symmetry of the binomial coefficients, the first and third sums
are both $V_n(u)$. Since $n=2u+t$, the middle sum contains exactly
$t-1$ terms and is nonempty for $t\ge3$. Hence
\[
 2^n>2V_n(u),
 \qquad\text{so}\qquad
 V_n(u)<2^{n-1}=M_t(n).
\]
Thus, in this case, the maximum size of a code avoiding the single
distance $n-t+1$ is strictly larger than the maximum size of a
binary $t$-intersecting code.

Now suppose that $n-t+1$ is even, and put $N=n-1$ and
$r=(n-t-1)/2$, so that $n-t=2r+1$. For any
$c\in Q_N$, the code $B_r(c)\times\{0,1\}$ has diameter at most
$2r+1=n-t$, and hence is $t$-intersecting. Its size is
$2V_N(r)=M_t(n)$. Since every binary $t$-intersecting code avoids
distance $n-t+1$,
its size is at most $M_t(n)$.
Thus, in the even case, the maximum size of a binary
$t$-intersecting code is exactly $M_t(n)$. We next show that the
extremal structures of the two problems are nevertheless different.

Indeed, not every maximum avoiding code is $t$-intersecting. In
\eqref{eq:extremal}, take $c_0=\mathbf0$ and $c_1=\mathbf1$. The
resulting maximum code contains $\iota_0(\mathbf0)$ and
$\iota_1(\mathbf1)$. Their first $N$ coordinates differ in every
position. If $N$ is even, their appended coordinates are also
different, so their Hamming distance is $N+1=n$; if $N$ is odd,
their appended coordinates are equal, so their Hamming distance is
$N=n-1$. In either case, this distance exceeds $n-t$, and hence the
code is not $t$-intersecting.

\par\medskip
\noindent\textbf{Funding.} This work was supported by the National Natural Science Foundation of China under Grant No.~12071351.

\noindent\textbf{Declaration of competing interests.} The authors declare no competing interests.

\section*{Declaration of generative AI and AI-assisted technologies in the manuscript preparation process}
During the early preparation of this manuscript, generative AI was used to assist with drafting, language refinement, and exploratory suggestions. The latest available DeepSeek model was used in the final stage to assist in checking the logical correctness of proofs of selected auxiliary lemmas and polishing the manuscript’s language. The central mathematical ideas, main results, and overall proof strategy were developed by the authors. All AI-assisted content was reviewed and verified by the authors, who take full responsibility for the manuscript.

\end{document}